\documentclass[letterpaper,12pt,oneside]{article}

\usepackage{amssymb}
\usepackage{amssymb}
\usepackage{enumerate}
\usepackage{amsmath}
\usepackage{amsthm}
\usepackage{textcomp}
\usepackage{graphicx,epsfig,tikz,mdwlist}
\usepackage{amsmath,amsfonts,amssymb,enumerate}
\usepackage{latexsym,amsfonts,amssymb,amscd}
\usepackage{fancybox,shadow,color,graphicx,psfrag,float}
\newcommand{\m}{\mathcal}

\newcommand{\F}{\mathbb{F}}

\newtheorem{thm}{Theorem}[section]
\newtheorem{pro}[thm]{Proposition}

\newtheorem{rem}[thm]{Remark}

\theoremstyle{definition}
\newtheorem{example}[thm]{Example}

\usepackage[colorlinks=true,linkcolor=red,
citecolor=red, urlcolor=blue]{hyperref}

\title{A note on subtowers and supertowers of recursive towers of function fields\footnote{This work was partially supported by CONICET and UNL CAI+D 2016.}}

\author{M. Chara\footnote{Consejo Nacional de Investigaciones Cient\'ificas y T\'ecnicas,  Argentina;  Universidad Nacional del Litoral, Santa Fe, Argentina; mchara@santafe-conicet.gov.ar}; H. Navarro\footnote{Universidad del Valle, Cali, Colombia; horacio.navarro@correounivalle.edu.co} and 
R.Toledano\footnote{Universidad Nacional del Litoral, Santa Fe, Argentina; ridatole@gmail.com}}

\date{\footnotesize{ }}

\begin{document}
\maketitle

\begin{abstract}

In this paper we  study the problem of constructing non-trivial subtowers and supertowers of recursive towers of function fields over finite fields.
\end{abstract}

\section{Introduction}

Let $q$ be a prime power and let $F$ be an algebraic function field of one variable over a finite field $\F_q$ of cardinality $q$. In \cite{ih81}, Ihara introduced the function $$A(q)=\limsup_{g \rightarrow \infty} \frac{N_q(g)}{g},$$
where $N_q(g)$ is the maximum number of rational places that a function field over $\F_q$ of genus $g$ can have.  This function measures how large  the number of rational places in function fields with respect to their genus can be and it shows up, for instance,  in the so called Tsfasman, Vladut and Zink bound in coding theory (see, for example, \cite[Proposition 8.4.6]{Stichbook09}). This is a good motivation to find the exact values of $A(q)$ but nothing is known except that $A(q^2)=q-1$. 

Because finding the values of $A(q)$ when $q$ is not a square has proven to be really hard,  most efforts are directed to give lower bounds for  $A(q)$.    One way of obtaining non-trivial lower bounds for $A(q)$ is through the construction of asymptotically good towers of function fields over $\F_q$. Following \cite{Stichbook09} a tower is a sequence $\m{F}=(F_0, F_1, \ldots)$ of function fields over a fixed finite field $\F_q$, such that for each $n\geq 0$ the extension  $F_{n+1}/F_n$ is  finite and separable,  $\F_q$ is the full constant field of $F_n$ and the genus $g(F_n)$ of $F_n$ goes to infinity along with $n$. If $N(F_n)$ denotes the number of $\F_q$-rational places of $F_n$, then the limit $\lambda(\m{F})=\lim_{n\rightarrow\infty}N(F_n)/g(F_n)$ exists and it is called the limit of the tower. Clearly, this limit provides a lower bound for the quantity $A(q)$ and when $\lambda(\m{F})=A(q)$ the tower $\m{F}$ is called optimal over $\F_q$.

By using only basic results of valuation theory and ramification in Artin-Schreier extensions, Garcia and Stichtenoth gave in \cite{GS95} the first example of an optimal recursive tower of function fields over $\F_{q^2}$. Recursive means that all the extensions are defined by the same equation (see Section \ref{basics} for details) and the interest in finding good recursive towers lies in the possibility of having a concrete description of the geometric Goppa codes attached to them (see \cite{Stichbook09}).

 One tricky thing when working with these recursive towers is that many times  different equations give rise to the same tower and it is not trivial at all
  how to decide if the chosen equation is the best one to work with, in the sense that this equation may not be the most suitable for the determination of some invariants in the tower. With this in mind, the concepts of subtowers and supertowers gain importance. Basically, a subtower $\m{E}=(E_0, E_1, \ldots)$ of a tower $\m{F}=(F_0, F_1, \ldots)$ is a tower in which each function field $E_i$ is embedded in some $F_j$, for $j\geq i$. Equivalently it is also said that $\m{F}$ is a supertower of $\m{E}$ and we always have that $\lambda(\m{E})\geq \lambda(\m{F})$. (See Section~\ref{secsubsequences} for precise definitions). In this regard it is important to recognize when two equations define the same tower and also if an equation defines a supertower or a subtower of an already studied tower or an easier one to study. {The above definitions can also be given in the slightly weaker situation of sequences  $\m{F}=(F_0, F_1, \ldots)$  of function fields (see Section \ref{basics}) where the condition  $g(F_n)\rightarrow \infty$  as $n\rightarrow\infty$ is not required.}

 {The aim of this paper is to provide a systematic method to construct recursive subsequences and supersequences of function fields from recursive sequences and to check if two apparently different equations give rise to recursive sequences where one is a subsequence of the other. This is done in Section~\ref{secsubsequences} and in Theorem~\ref{method} we prove that our method actually gives rise to a proper recursive subsequence of a given recursive sequence. An interesting feature of these results is that they can be easily implemented in a computer so we were able to search for many equations defining subsequences.}

  The paper is organized as follows: in Section~\ref{basics} we give some basic definitions.  In Section~\ref{secsubsequences} we present our main results, which will be used in Section~\ref{examples} to work with different examples using our method.  Finally we show in Section~\ref{sec5} an interesting application of our results by finding an optimal quadratic recursive tower $\m{E}=(E_0,E_1,\ldots)$ over $\F_4$ whose field extensions $E_{i+1}/E_i$ are Artin-Schreier extensions but the tower $\m{E}$ itself is not recursively defined by an Artin-Schreier equation.

  \section{Basic definitions}\label{basics}
  Following \cite{GS07} and \cite{Stichbook09} by a recursive sequence of function fields over $\F_q$ we mean  a sequence of function fields $\m{F}=(F_0, F_1, \ldots)$ over $\F_q$,  a sequence $\{x_i\}_{i=0}^{\infty}$ of transcendental elements over $\F_q$ and a bivariate polynomial
\[H\in\F_q[S,T]\,,\]
 such that
\begin{enumerate}
\item $F_0=\F_q(x_0)$,
\item $F_{i+1}=F_i(x_{i+1})$ where $H(x_i,x_{i+1})=0$ for $i\geq 0$, and
\item the polynomial $H(x_i,T)\in F_i[T]$ is separable for $i\geq 0$.
\end{enumerate}
Notice that from this definition we have that each field extension $F_{i+1}/F_i$ is finite (because $[F_{i+1}:F_i]\leq \deg_T(H(x_i,T))$) and separable. Also \[F_i=\F_q(x_0,\ldots,x_i)\quad\text{for $i\geq 0$}\,,\] so that
\[F_0=\F_q(x_0)\subseteq F_1\subseteq\ldots F_i\subseteq F_{i+1}\subseteq\ldots\]

We shall say that a recursive sequence of function fields is non-trivial if $[F_{i+1}:F_i]\geq 2$ for every $i\geq 0$, in other words $F_i\subsetneq F_{i+1}$.

If $\m{F}$ is a non-trivial recursive sequence such that the genus $g(F_i)\rightarrow \infty$ as $i\rightarrow\infty$ and $\F_q$ is algebraically closed in each $F_i$ we shall say that $\m{F}=(F_0, F_1, \ldots)$ is a recursive tower of function fields over $\F_q$. As stated in \cite{Stichbook09}, it suffices to have that $g(F_i)\geq 2$ for some index $i\geq 1$ in order to have that $g(F_i)\rightarrow \infty$ as $i\rightarrow \infty$ . When $\F_q$ is algebraically closed in each $F_i$ it is customary to say that $\F_q$ is the full field of constants of each $F_i$.

The following definitions are important when dealing with the asymptotic behaviour of a tower. Let $\m{F}=(F_0, F_1, \ldots)$ be a  tower (not necessarily recursive) of function fields
over a finite field $\F_q$. Let $N(F_i)$ be the number of rational places of $F_i$. The  splitting rate $\nu(\m{F})$ and the genus $\gamma(\m{F})$ of $\m{F}$ over $F_0$ are defined, respectively, as
$$\nu(\m{F})\colon=\lim_{i\rightarrow \infty}\frac{N(F_i)}{[F_i:F_0]}\,, \qquad\gamma(\m{F})\colon=\lim_{i\rightarrow \infty}\frac{g(F_i)}{[F_i:F_0]}\,.$$
The  limit $\lambda(\m{F})$ of  the tower $\m{F}$ is defined as $$\lambda(\m{F})\colon=\lim_{i\rightarrow \infty}\frac{N(F_i)}{g(F_i)}\,.$$
It can be seen that all the above limits exist and that $A(q)\geq\lambda(\m{F})\geq 0$ (see \cite[Chapter 7]{Stichbook09}). The tower $\m{F}$ is called asymptotically good (over $\F_q$) if $\lambda(\m{F})>0$ (in particular $\m{F}$ is called asymptotically optimal over $\F_q$ if $\lambda(\m{F})=A(q)$). Otherwise  $\m{F}$ is called asymptotically bad.

 If a tower $\m{F}=(F_0, F_1, \ldots)$ of function fields over $\F_q$ is recursively defined by a polynomial of the form
\begin{equation}\label{abtower}
H(S,T):=a_1(T)b_2(S)-a_2(T)b_1(S),
\end{equation}
where $a_1$, $a_2$, $b_1$, $b_2 \in \F_q[T]$ are polynomials such that
\[\gcd(a_1,a_2)=\gcd(b_1,b_2)=1,\]
and
\[[F_{i+1}:F_i]=\deg_TH, \]
 we shall say that $\mathcal{F}$ is an $(a,b)$-recursive tower of function fields over $\F_q$ in order to make reference to the rational functions
\[a(T):=\frac{a_1(T)}{a_2(T)}\quad\text{and}\quad b(S):=\frac{b_1(S)}{b_2(S)}\,,\]
defining the sequence. A tower recursively defined by an equation with mixed variables is a recursive tower which is not an $(a,b)$-tower.

Of course not any choice of rational functions $a,b\in\F_q(T)$ will give rise to a recursive tower over $\F_q$. For example, it was shown in \cite{MaWu05} that absolutely irreducible and symmetric polynomials\footnote{This  means that $H$ is irreducible over an algebraic closure of $\F_q$ and that $H(S,T)=H(T,S).$} $H\in \F_q[S,T]$  do not give rise to towers if the extension $\F_q(x,y)/\F_q(x)$ is Galois where $H(x,y)=0$ and $x$ is transcendental over $\F_q$. They actually proved that, under the above conditions, if $F_{i+1}=F_i(x_{i+1})$ where $H(x_i,x_{i+1})=0$ for $i\geq 0$ then $F_i\subseteq F_1$ for all $i\geq 1$.

We now introduce a special family of equations defining sequences of function fields over a finite field $\F_q$ which will be used in the last section. 

\begin{pro}\label{propo3.31tesis}
Let $\m{F}=(F_0,F_1,\ldots)$ be a recursive sequence of function fields defined over a finite field $\F_q$ by the equation
\begin{equation}\label{eq1prop2.3}
y^m+f(x)y^{m-1}+\cdots+f(x)^{m-1}y+h(x)=0 
\end{equation}
where $m$ is a power of the characteristic $p$ of $\F_q$, $f$ and $h$ are the following rational functions 
$$f(T)=\frac{T-\gamma}{\alpha T-\beta}, \qquad h(T)=\frac{(T-\gamma)^j}{h_1(T)}-\gamma,$$
 with $h_1(T) \in \F_q[T]$, $\alpha, \beta, \in \F_q$ and $\gamma \in \F_m^{*}$, $j$ is a fixed positive integer such that $1<j\leq m$, and $T-\gamma$ is coprime with $(\alpha T-\beta) h_1(T).$ Then the simple zero $P_\gamma$ of $x_0-\gamma$ in  $F_0$ is totally ramified in the sequence. In particular all the extensions $F_{i+1}/F_i$ have degree $m$ and $\F_q$ is the full constant field of $F_i$ for every $i\geq0$.
\end{pro}
\begin{proof}
Notice that $y\neq f(x)$ in \eqref{eq1prop2.3} so that by multiplying \eqref{eq1prop2.3} by $y-f(x)$ we can rewrite \eqref{eq1prop2.3} as	
 \begin{equation}\label{eq2prop2.3}
 y^{m+1}+(h(x)-f(x)^m)y=h(x)f(x)
 \end{equation}

Thus we have a sequence $\{x_i\}_{i=0}^{\infty}$ of transcendental elements over $\F_q$ such that $F_{i+1}=F_i(x_{i+1})$ where $x_{i+1}$ is a root of	
 $$\phi(T)=T^{m+1}+(h(x_i)-f(x_i)^m)T-h(x_i)f(x_i)\in F_i[T],$$
 where $x_{i+1}\neq f(x_i)$ for all $i\geq 0$.
  
Assume now that $P$ is a simple zero of $x_i-\gamma$ in $F_i$ and let $Q$ be a place of $F_{i+1}$ lying over $P$. We will prove that $Q|P$ is totally ramified in $F_{i+1}/F_i$. For simplicity we write $y=x_{i+1}$ and $x=x_i$. Then the extension $F_{i+1}/F_i$ is defined by the equation $$y^{m+1}+(h(x)-f(x)^m)y=h(x)f(x),$$ and we claim that 
\begin{enumerate}[(i)]
	\item $\nu_{P}(h(x)f(x))=1$, 
	\item $\nu_{P}(h(x)+\gamma)=j$ and
\item $\nu_{P}(f(x))=1$. 	
\end{enumerate}
To see this notice that from our hypothesis $h_1(T)=(T-\gamma)h_2(T)+c$ where $h_2\in \F_q[T]$ and $c\in \F_q^*$. Then $\nu_P(h_1(x))=0$ and (ii) follows. We also have that $\nu_P(h(x))=0$. Now $\nu_{P}(f(x))=1$ because $P$ is a simple zero of $x-\gamma$ and if $\alpha\neq 0$ we write $\alpha x-\beta=\alpha(x-\gamma+\gamma-\beta/\alpha)$ with $\gamma-\beta/\alpha\neq 0$ by hypothesis. Thus (ii) and (iii) follow.

Therefore from \eqref{eq2prop2.3} we deduce that $$\nu_Q(y)+\nu_Q(y^m+h(x)-f(x)^m)=e(Q|P).$$
Let us see now that $\nu_Q(y)=0$: if $\nu_Q(y)<0$ then we have that $$\nu_Q(y^m+h(x)-f(x)^m)=m\nu_Q(y),$$
because $\nu_Q(h(x)-f(x)^m)=0$ by (iii) and (ii) above and this implies that 
$$(m+1)\nu_Q(y)=\nu_Q(y)+\nu_Q(y^m+h(x)-f(x)^m)=e(Q|P),$$ which is a contradiction. Similarly if  $\nu_Q(y)>0$ then we have that   $$\nu_Q(f(x)^iy^{m-i})=i\nu_Q(f(x))+(m-i)\nu_Q(y)>0,$$
 for each $0\leq i \leq m-1$. Now from \eqref{eq1prop2.3}, since $\nu_Q(h(x))=0$, we see  that  
$$m\nu_Q(y)=\nu_Q(-f(x)y^{m-1}-f(x)^2y^{m-2}-\cdots-f(x)^{m-1}y-h(x))=0,$$ which is again a contradiction. 

Thus $\nu_Q(y)=0$ and then $\nu_Q(f(x)^iy^{m-i})=ie(Q|P)$ for $1\leq i \leq m$.  Also, since $\gamma\in \F_m^*$ we have
$$(y-\gamma)^m=-f(x)y^{m-1}-f(x)^2y^{m-2}-\cdots-f(x)^{m-1}y-(h(x)+\gamma)$$ so that 
$$m\nu_Q(y-\gamma)=e(Q|P).$$ Therefore $Q|P$ is totally ramified in $F_{i+1}/F_i$ and we also see that $Q$ is a simple zero of $y-\gamma$ in $F_{i+1}$.  Since $P_{\gamma}$ is a simple zero of $x_0-\gamma$ in the rational function field $F_0$, the result follows from an  inductive argument. \end{proof}

\section{Constructing subsequences and supersequences}\label{secsubsequences}

Let $\m{F}=(F_0, F_1, \ldots)$ be a sequence of function fields over $\F_q$. A sequence $\mathcal{E}=(E_0, E_1, \ldots)$ of function fields over $\F_q$ is called {\em subsequence} if for each $i\geq0$ there exists an index $j=j(i)$ and an embedding $\varphi_i:E_i\rightarrow F_j$ over $\F_q$. If, in addition, $\varphi_i(E_i)\subsetneq F_j$ for infinitely many  $i\geq 0$ we shall say that $\mathcal{E}$ is a proper subsequence of $\mathcal{F}$.  Moreover if the sequences $\m{F}=(F_0, F_1, \ldots)$ and $\mathcal{E}=(E_0, E_1, \ldots)$ are actually  towers of function fields it is said that either $\mathcal{E}$ is a {\em subtower} of $\mathcal{F}$ or that $\mathcal{F}$ is a {\em supertower} of $\mathcal{E}$.
\vspace{0.3cm}

From now on we will always assume  that a rational function $a\in \F_q(T)$ is given in its lowest terms, i.e.  there are two coprime polynomials $a_1$, $a_2\in \F_q[T]$ such that $a=a_1/a_2$. In this case we define {\em the degree} of  $a\in \F_q(T)$  as $\deg(a):=\max\{\deg(a_1),\deg(a_2)\}$.
\vspace{0.3cm}

We prove next a preliminary result which will be important for our method to construct a recursive subsequence from a given $(a,b)$-recursive sequence $\m{F}$.

\begin{pro}\label{proposub1}
Let $\m{F}=(F_0, F_1, \ldots)$ and $\m{E}=(E_0, E_1, \ldots)$ be non-trivial  recursive sequences of function fields over $\F_q$ defined respectively by the equations 
$$f(x_i, x_{i+1})=0\qquad \text{and}\qquad h(y_i, y_{i+1})=0$$ where $f(X,Y)$ and $h(X,Y)$ are bivariate polynomials with coefficients in $\F_q$ and $\{x_i\}_{i\geq 0}$ and $\{y_i\}_{i\geq 0}$ are sequences of transcendental elements over $\F_q$. For each $i\geq 0$, let us assume that $$\deg_Y h \leq [F_{i+1}:F_i]$$ and that $y_i=g(x_i)$ with $g(T)=g_1(T)/g_2(T) \in \F_q(T)$ a  rational function of degree bigger than one. Then $E_i \subsetneq F_i$ for every $i\geq 0$, i.e., $\m{E}$ is a proper subsequence of $\m{F}$. 
\end{pro}
\begin{proof}
Since $y_i=g(x_i)$, $E_i=\F_q(y_0,\dots, y_i)$ and $F_i=\F_q(x_0, \ldots, x_i)$ we have that $E_i \subseteq F_i$ for every $i\geq 0$. On the other hand, we can assume without loss of generality, that the degree of the rational function $g(T)$ is $\deg g=\deg g_1$. For each index $i\geq 0$, let us consider the polynomial 
$$\varphi_i(T)=g_1(T)-g_2(T)y_i \in E_i[T].$$

It is clear that the element $x_i$ is a root of the polynomial $\varphi_i(T)$ for each $i\geq 0$. Since $E_0=\F_q(y_0)$ and $F_0=\F_q(x_0)$ then (see Section 14.9 of \cite{DF}) $[F_0:E_0]=\deg g>1$. Now let $d_i=[F_i:E_i]$. Then we have 
$$d_{i+1}[E_{i+1}:E_i]=d_i[F_{i+1}:F_i].$$ We will show by induction that $d_i>1$ for $i\geq 1$. If $d_1=1$, then $$[E_1:E_0]=d_0 [F_1:F_0] =\deg g [F_1:F_0].$$ By hypothesis, we have that $[E_1:E_0]\leq \deg_Y h \leq [F_1:F_0]$ and therefore $$\deg g [F_1:F_0] = [E_1:E_0] \leq [F_1:F_0],$$
and so $\deg g \leq 1$ which is a contradiction. Suppose now that $d_i>1$ and that $d_{i+1}=1$. Then, 
$$d_i[F_{i+1}:F_i]=[E_{i+1}:E_i]\leq \deg_Y h\leq [F_{i+1}:F_i]$$ and thus $d_i\leq 1$ which is a contradiction. Therefore $E_i \subsetneq F_i$ for each $i~\geq~0$, i.e., $\m{E}$ is a proper subsequence of $\m{F}$.  \end{proof}

\begin{rem}
From the proof of the previous theorem, we have that $$\deg g=[F_0:E_0]\qquad\text{and}\qquad [F_i:E_i]\geq \deg g,$$ for each $i\geq 1$. Moreover, equality  $[F_i:E_i]=\deg g$ holds for every $i\geq 0$ if the equality $[E_{i+1}:E_i]=\deg_Y h=[F_{i+1}:F_i]$ is assumed. 
\end{rem}

Now we are in a position to prove the main result of this section.
\begin{thm}[The method] \label{method}
Let $\m{F}=(F_0, F_1, \ldots)$ be a non-trivial (a,b)-sequence of function fields over $\F_q$. Let us assume that $A(T)$, $B(T)$, $g(T)$, $s(T) \in \F_q(T)$ are rational functions such that
\begin{equation}\label{funequsubtower}
A\circ g= s\circ a \qquad \text{and} \qquad B\circ g=s \circ b.
\end{equation} Then the sequence $\m{E}=(E_0, E_1, \ldots)$ recursively defined by $E_0=\F_q(g(x_0))$ and $E_{i+1}=E_i(g(x_{i+1}))$, with $A(g(x_{i+1}))=B(g(x_i))$, is a recursive (A,B)-subsequence of $\m{F}$. If $\m{E}$ is a non-trivial sequence and for every $i\geq 0$ we have that $$\deg A\leq [F_{i+1}:F_i]$$ then $\m{E}$ is a proper subsequence of $\m{F}$.
\end{thm}

\begin{proof}
Let $\{x_i\}_{i=0}^{\infty}$ be a sequence of transcendental elements over $\F_q$ such that $a(x_{i+1})=b(x_i)$ for each $i\geq 0$. Then, $\{g(x_i)\}_{i=0}^{\infty}$ is a sequence of transcendental elements over $\F_q$  and we have that 
$$A(g(x_{i+1}))=s(a(x_{i+1}))=s(b(x_i))=B(g(x_i)).$$
It is clear that $E_i\subseteq F_i$ for each $i\geq 0$ so that $\m{E}$ is a subsequence of $\m{F}$. 

Now if $\m{E}$ is non-trivial and $\deg A\leq [F_{i+1}:F_i]$ then, for every $i\geq 0$, $$[E_{i+1}:E_i]\leq \deg A\leq [F_{i+1}:F_i],$$ and by Proposition \ref{proposub1} we have that $\m{E}$ is a proper subsequence of $\m{F}$ as desired. 
\end{proof}

\begin{rem}
From the last part of the proof of Theorem \ref{method} we have that \[[F_i:E_i]=\deg g,\] for every $i\geq 0$ if $\deg A=[F_{i+1}:F_i]$ for every $i\geq 0$.
\end{rem}

\begin{rem}
Notice that condition $\deg A\leq [F_{i+1}:F_i]$ in Theorem~\ref{method} can be replaced by $[E_{i+1}:E_i]\leq [F_{i+1}:F_i]$ if $\deg A> [F_{i+1}:F_i]$.
\end{rem}

\begin{rem}
If $\m{F}$ is an $(a,b)$-tower and there is a $(A,B)$-subtower $\m{E}$ and functions $g$ and $s$ such that the conditions in Theorem \ref{method} hold, then for any rational function $$f(T)=\frac{aT+b}{cT+d}\qquad \text{with }a,b,c,d \in \F_q \quad \text{and}\quad ad-cb\neq 0$$ we also have that the functions $A\circ f$, $B\circ f$, $f^{-1}\circ g$ and $s$ satisfy the conditions in Theorem \ref{method}.  Therefore 
if we define $$G_0=\F_q(f^{-1}(g(x_0))) \qquad  \text{and}\qquad   G_{i+1}=G_{i}(f^{-1}(g(x_{i+1})))\quad \text{for }i\geq 0$$ then we obtain an $(A\circ f,B\circ f)$-subtower $\mathcal{G}=(G_0, G_1, \ldots)$ of function fields over $\F_q$ of $\mathcal{F}$. Actually, $\m{E}$ and $\m{G}$ are the same tower (See \cite[Equation (2.3)]{BGS06}). 

\end{rem}

\section{Examples}\label{examples} We will show next that many subtowers studied  in the literature can be obtained using our method presented in Theorem \ref{method}.
First we list some well known recursive towers.

\begin{enumerate}[1)]
\item(Bezerra and Garcia \cite{BG04}) The equation 
\begin{equation}\label{ec.0}
\frac{y-1}{y^q}=\frac{x^q-1}{x}
\end{equation}
defines an asymptotically optimal tower $\m{F}_0$ over the field $\F_{q^2}$.

\item (Garcia, Stichtenoth \cite{GS00}) The equation 
\begin{equation}\label{ec.1}
y^q+y=\frac{x^q}{x^{q-1}+1}
\end{equation}
defines an asymptotically optimal tower $\m{F}_1$ over the field $\F_{q^2}$.

\item(Bassa et al. \cite{BaGS07}) The equation 
\begin{equation}\label{ec.2}
(y^q-y)^{q-1}+1=-\frac{x^{q(q-1)}}{(x^{q-1}-1)^{q-1}}\,,
\end{equation}
defines an asymptotically good tower $\m{F}_2$ over the field $\F_{q^3}$.

\item(Bezerra et al. \cite{BGS03}) The equation  
\begin{equation}\label{ec.3}
\frac{1-y}{y^q}=\frac{x^q+x-1}{x}\,,
\end{equation}
defines an asymptotically good tower $\m{F}_3$ over the field $\F_{q^3}$.

\item (Caro, Garcia \cite{carogarcia}) The equation
\begin{equation}\label{ec.4} 
y^{q+1}+y=\frac{x+1}{x^{q+1}},
\end{equation}
defines an asymptotically good tower $\m{F}_4$ over the field $\F_{q^3}$.

\item(Garcia et al. \cite{GSR03}) Let $q=p^{2n}$ where $p$ is an odd prime. The equation of Kummer type
\begin{equation}
y^2=\frac{x^2+1}{2x}
\end{equation}
defines an asymptotically good tower  $\m{F}_5$ over $\F_q$.

\item(Garcia et al. \cite{GSR03}) The equation of Kummer type  $$y^2=\frac{x^2}{x-1},$$
defines an asymptotically good tower  $\m{F}_6$ over $\F_9$.
\end{enumerate}

\begin{example}
In 2004 Bezerra and Garcia proved in \cite{BG04} that the tower $\m{F}_0=(E_0, E_1,\ldots)$  is a subtower of the tower $\m{F}_1=(F_0,F_1,\ldots)$.

 Actually, the $(A,B)$-tower $\m{F}_0$ is a proper subtower of the $(a,b)$-tower $\m{F}_1$ as can be seen from Theorem~\ref{method}. Let us consider $$g(T)=\frac{1}{T^{q-1}+1}\qquad \text{and}\qquad s(T)=-T^{q-1}.$$ Since $$a(T)=T^q+T, \quad b(T)=-\frac{T^q}{T^{q-1}+1}, \quad A(T)=\frac{T-1}{T^q}\quad \text{and}\quad B(T)=\frac{T^q-1}{T}, $$
then it is easy to check that $$A\circ g= s\circ a \qquad \text{and} \qquad B\circ g=s \circ b.$$
Moreover, $\deg A=[E_{i+1}:E_i]=q=[F_{i+1}:F_i]=\deg a$ and therefore $[F_i:E_i]=\deg g=q-1$ and $\m{F}_0$ is a proper subtower of $\m{F}_1$.
\end{example}

\begin{example} In this example we show  that $\m{F}_3$ is a subtower of  $\m{F}_2$. The tower $\m{F}_2=(F_0,F_1,\ldots)$ is an $(a,b)$-tower where  
$$a(T)=(T^q-T)^{q-1}+1 \quad \text{ and } \quad b(T)=\frac{T^{q(q-1)}}{(T^{q-1}-1)^{q-1}}$$
and satisfy  $[F_1:F_0]=q(q-1)$ and $[F_{i+1}:F_i]=q$ for all $i\geq 1$. The tower $\m{F}_3=(H_0, H_1, \ldots)$ is an $(A,B)$-tower with
$$ A(T)=\frac{1-T}{T^{q}} \quad \text{ and } \quad {B}(T)=\frac{T^{q}+T-1}{T}\cdot$$   
In this case $[H_{i+1}:H_i]= q$ for all $i\geq 0$. Let us consider  $g(T)=-\frac{1}{T^{q-1}-1}$ and $s(T)=-T+1$ then
 $$A\circ g= s\circ a \qquad \text{and} \qquad B\circ g=s \circ b,$$ and $\deg A=[H_{i+1}:H_i]=q$.
From Theorem~\ref{method} we have that   $\m{F}_3$ is  a proper subtower of $\mathcal{F}_2$ over $\F_{q^3}$.

Notice that in this case $$[F_0:H_0]=q-1 \quad \text{and}\quad [F_i:H_i]=(q-1)^2\quad \text{for }i\geq 1.$$ Thus, each $F_{i+1}$ is actually the composition field of the fields $F_i$ and $H_{i+1}$.
\end{example}

\begin{example}In this example we show that the $(A,B)$-tower $\m{F}_4=(G_0, G_1, \ldots)$ where
$$A(T)=T^{q+1}+T \quad \text{ and } \quad B(T)=\frac{T+1}{T^{q+1}}\cdot$$
and $[G_{i+1}:G_i]=q$ for all $i\geq 0$ is a subtower of  $\m{F}_3=(H_0, H_1, \ldots)$ of the previous example. It is not hard to check that conditions in Theorem \ref{method} hold if we choose $$g(T)=\frac{T^q+T-1}{1-T}\qquad \text{and}\qquad s(T)=\frac{1-T}{T}+\left(\frac{1-T}{T}\right)^{q+1}$$.

\end{example}

\begin{example} Let $q=p^{2n}$ where $p$ is an odd prime. The equation of Kummer type
\begin{equation}
y^2=\frac{x^2+1}{2x}\,,
\end{equation}
defines the $(a,b)$-recursive tower $\m{F}_5$ of function fields over $\F_q$. In this case
\[F_{i+i}=F_i(x_{i+1})\quad\text{with}\quad x_{i+1}^2=\frac{x_i^2+1}{2x_i}\quad\text{for $i\geq 0$}\,,\]
and we have that $a(T)=T^2$ and $b(T)=(T^2+1)/2T$.

Now if $g(T)=2T^2$, $s(T)=4T^4$, $A(T)=T^2$ and $B(T)=(T+2)^2/2T$ then it is easy to check that \[A\circ g= s\circ a \qquad \text{and} \qquad B\circ g=s \circ b,\]
so that the equation
\[y^2=\frac{(x+2)^2}{2x}\,,\]
defines an $(A,B)$-recursive proper subsequence $\mathcal{E}=(E_0, E_1, \ldots)$ of $\mathcal{F}_5$ over $\F_q$.  
 In fact, $\mathcal{E}$ is actually a proper subtower of $\mathcal{F}_5$ over $\F_q$.  This subtower was also obtained in \cite{MaWu05} using a method due to Elkies.
\end{example}

\begin{example}
Now we want to determinate whether the tower $\m{F}_6$ over $\F_9$ recursively defined by $$y^2=\frac{x^2}{x-1},$$
has any relationship with some of the already known asymptotically good towers over $\F_9$. 
We perform a computational search of possible functions $g(T)$ and $s(T)$ described in our method of Theorem \ref{method} with $A(T)=T^2$, $B(T)=T^2/(T-1)$ and some known $({a},{b})$-tower
over $\F_9$. As a result we have that using $${a}(T)=T^2,\quad {b}(T)=\frac{(T+2)^2}{2T},\quad g(T)=T^2+1 \quad \text{and}\quad  s(T)=(T+1)^2\, ,$$ equation \eqref{funequsubtower} is satisfied and also Theorem~\ref{method} holds.

Therefore the tower $\m{F}_6$ is actually a subtower of the tower $\m{E}$ in the previous example and therefore is also a subtower of $\m{F}_5$.

 Notice that the tower $\mathcal{F}_6$ was studied in \cite{GSR03} but it was not mentioned that $\mathcal{F}_6$ is a subtower of $\mathcal{E}$  and $\mathcal{F}_5$ over $\F_9$. Moreover, performing the change of variables $x_1=1/x$ and $y_1=1/y$ we get the Fermat type tower recursively defined by $$y_1^2=x_1(1-x_1).$$ Therefore this example was not new.
\end{example}

\section{An optimal recursive quadratic tower  over $\F_4$ with mixed variables.}\label{sec5}

Let us consider now the additive polynomial $\phi=T^2+(x+1)T$ over the rational function field $\F_4(x)$. The aim of this section is to show in a simple way that the polynomial $\phi$ gives rise to an asymptotically optimal tower $\m{E}$ over $\F_4$ whose function field extensions are Artin-Schreier extensions. We will prove it by showing that the tower $\m{F}_1$ over $\F_4$ presented in the previous section is a supertower for $\m{E}$ using Theorem \eqref{method}. 

Let  $\m{E}=(E_0, E_1, \ldots)$ be  the sequence recursively defined over $\F_4$ by the equation with mixed variables $\phi(y)=x^2$, i.e.
\begin{equation}\label{eq1sec5}
y^2+(x+1)y=x^2.
\end{equation} 
 This equation defines a non-trivial sequence of function fields because we are in the hypothesis of Proposition~\ref{propo3.31tesis} with $\gamma=1$, $\alpha=0$, $\beta=1$ and the polynomials 
\[f(T)=T+1\qquad \text{and}\qquad h(T)=T^2=(T+1)^2+1.\] 
Therefore the simple zero $P_1$ of $x_0+1$ in the rational function field $E_0=\F_4(x_0)$ is totally ramified in $\m{E}$ so that  $\F_4$ is the full field of constant of each field $E_i$ and  every extension $E_{i+1}/E_i$ has degree $2$. 

Notice that each extension $E_{i+1}/E_i$ is an Artin-Schreier extension defined by the equation
\begin{equation}\label{eq2sec5}
z^2+z=\frac{x^2_i}{x^2_i+1},
\end{equation}
where $z=x_{i+1}/(x_i+1)$. However the sequence $\m{E}$ is not recursively defined by \eqref{eq2sec5} because $z$ is in terms of $x_i$ and $x_{i+1}$.   From Kummer's theorem and the theory of Artin-Schreier extensions (see Chapter 3 of \cite{Stichbook09}) it is easily seen that $P_1$ is the only place of $E_0$ ramified in $E_1$, the zero $P_0$  and the pole $P_{\infty}$ of $x_0$ in $E_0$ split completely  into a simple zero of $x_{1}$ and a simple zero of $x_{1}+1$ in $E_{1}$ and the genus of $E_1=\F_4(x_0,x_1) $ is zero. By applying repeatedly Kummer's theorem, we see from \eqref{eq1sec5} that for $i\geq 0$ there is a simple zero of $x_i$ in $E_i$ which splits into a simple zero of $x_{i+1}$ and a simple zero of $x_{i+1}+1$ in $E_{i+1}$. From the proof of Proposition~\ref{propo3.31tesis} we have that this simple zero of $x_{i+1}+1$ in $E_{i+1}$ is totally ramified in $E_j$ for $j>i+1$ and from the theory of Artin-Schreier extensions we have that the corresponding different exponents are 2.  With all of this and  Hurwitz's genus formula it is easy to verify that the genus $E_3$ is at least 3 so that $\m{E}$ is tower of function fields over $\F_4$ (it is also a tower over $\F_2$). 

Now we will prove that $\m{E}$ is an optimal tower over $\F_4$ by finding an optimal supertower using Theorem~\ref{method}. In this case, we have that $$A(T)=T^3+T\qquad \text{and}\qquad B(T)=T^3+T^2$$ and let us consider the functions
\[a(T)=T^2+T,\,\, b(T)=\frac{T^2}{T+1}, \,\, g(T)=\frac{T^2+T}{T^2+T+1}\,\, \text{ and }\,\, s(T)=\frac{T}{(T+1)^3}.\]
Then \[ A(g(T))=\frac{T^2+T}{T^6+T^5+T^3+T+1}=s(a(T)),\] and \[ B(g(T))=\frac{T^4+T^2}{T^6+T^5+T^3+T+1}=s(b(T)).\]  
so that the equation $$y^2+y=\frac{x^2}{x+1},$$ defines a recursive tower which is a supertower for $\m{E}$. But the above equation  is the optimal tower $\m{F}_1$   defined by \eqref{ec.1} over $\F_4$  and thus the optimality of $\m{E}$  over $\F_4$ follows.

\bibliographystyle{elsarticle-num}

\begin{thebibliography}{00}

\bibitem{BG04}
J.~Bezerra and A.~Garcia.
\newblock A tower with non-Galois steps which attains the Drinfeld-Vladut bound.
\newblock {\em Journal of Number Theory}, 106(1):142--154, 2004.


\bibitem{BaGS07}
A.~Bassa, A.~Garcia and H.~Stichtenoth.
\newblock  A New tower over cubic finite fields.
\newblock {\em Mosc. Math. J.}, 8(3):401--418, 2008.

\bibitem{BGS03}
J.~Bezerra, A.~Garcia and H.~Stichtenoth.
\newblock An explicit tower of function fields over cubic finite fields and Zink's lower bound.
\newblock {\em J. Reine Angew. Math.} 589:159--199, 2005.


\bibitem{BGS06}
P.~Beelen, A.~Garcia and H.~Stichtenoth.
\newblock Towards a classification of recursive towers of function fields over finite fields.
\newblock{\em Finite Fields Appl}, 12(1):56--77, 2006.


\bibitem{carogarcia}
N.~Caro and A.~Garcia.
\newblock On a tower of Ihara and its limit.
\newblock{\em Acta Arithmetica}, 151:191--200, 2012.

\bibitem{DF}
D.~Dummit and R.~Foote.
\newblock Abstract algebra.
\newblock Wiley. Third Ed. 2004. 

\bibitem{GS95}
A.~Garcia and H.~Stichtenoth.
\newblock A tower of Artin-Schreier extensions of function fields
attaining the Drinfeld-Vladut bound.
\newblock {\em Invent. Math.} 121, 211-222, 1995.

\bibitem{GS00}
A.~Garcia and H.~Stichtenoth.
\newblock On the asymptotic behaviour of some towers of function fields over finite fields. 
\newblock {\em Journal of Number Theory}, 61, 248--273, 1996.

\bibitem{GS07}
A.~Garcia and H.~Stichtenoth.
\newblock Explicit towers of function fields over finite fields.
\newblock In {\em Topics in geometry, coding theory and cryptography}, volume~6
  of {\em Algebr. Appl.}, pages 1--58. Springer, Dordrecht, 2007.

\bibitem{GSR03}
A.~Garcia, H.~Stichtenoth, and H.~R{\"u}ck.
\newblock On tame towers over finite fields.
\newblock {\em J. Reine Angew. Math.}, 557:53--80, 2003.

\bibitem{ih81}
Y.~Ihara.
\newblock Some remarks on the number of rational points of algebraic curves
  over finite fields.
\newblock {\em J. Fac. Sci. Univ. Tokyo Sect. IA Math.}, 28(3):721--724 (1982),
  1981.

\bibitem{MaWu05}
H.~Maharaj and J.~Wulftange.
\newblock On the construction of tame towers over finite fields.
\newblock {\em J. Pure Appl. Algebra}, 199(1-3):197--218, 2005.


\bibitem{Stichbook09}
H.~Stichtenoth.
\newblock Algebraic function fields and codes. 
\newblock GTM 254, Springer-Verlag, Berlin, 2nd Ed., 2009.


\end{thebibliography}

\end{document}